\documentclass[sn-mathphys-ay]{sn-jnl}

\usepackage{graphicx}%
\usepackage{multirow}%
\usepackage{amsmath,amssymb,amsfonts}%
\usepackage{amsthm}%
\usepackage{mathrsfs}%
\usepackage[title]{appendix}%
\usepackage{xcolor}%
\usepackage{textcomp}%
\usepackage{manyfoot}%
\usepackage{booktabs}%
\usepackage{algorithm}%
\usepackage{algorithmicx}%
\usepackage{algpseudocode}%
\usepackage{listings}%

\theoremstyle{thmstyleone}%
\newtheorem{theorem}{Theorem}
\newtheorem{lemma}{Lemma}
\newtheorem{proposition}[theorem]{Proposition}%
\newtheorem{corollary}[theorem]{Corollary}

\theoremstyle{thmstyletwo}%
\newtheorem{example}{Example}%
\newtheorem{remark}{Remark}%
\newtheorem{assume}{Assumption}

\theoremstyle{thmstylethree}%
\newtheorem{definition}{Definition}%

\begin{document}
	
	\title[Misclassification bounds for PAC-Bayesian sparse deep learning]{Misclassification bounds for PAC-Bayesian sparse deep learning}
	
	
	\author*{\fnm{The Tien} \sur{Mai}}\email{the.t.mai@ntnu.no}
		
	\affil{
		\orgdiv{Department of Mathematical Sciences}, 
		\orgname{Norwegian University of Science and Technology}, 
		\orgaddress{
			\city{Trondheim}, \postcode{7034}, 
			\country{Norway}}
		}

\abstract{
Recently, there has been a significant focus on exploring the theoretical aspects of deep learning, especially regarding its performance in classification tasks. Bayesian deep learning has emerged as a unified probabilistic framework, seeking to integrate deep learning with Bayesian methodologies seamlessly. However, there exists a gap in the theoretical understanding of Bayesian approaches in deep learning for classification. This study presents an attempt to bridge that gap. By leveraging PAC-Bayes bounds techniques, we present theoretical results on the prediction or misclassification error of a probabilistic approach utilizing Spike-and-Slab priors for sparse deep learning in classification. We establish non-asymptotic results for the prediction error. Additionally, we demonstrate that, by considering different architectures, our results can achieve minimax optimal rates in both low and high-dimensional settings, up to a logarithmic factor. Moreover, our additional logarithmic term yields slight improvements over previous works. Additionally, we propose and analyze an automated model selection approach aimed at optimally choosing a network architecture with guaranteed optimality.
}




\maketitle

\section{Introduction}

Deep learning has achieved remarkable success in various applications, including computer vision, natural language processing, and other areas related to pattern recognition and classification \citep{goodfellow2016deep,lecun2015deep}. As its applicability continues to expand, there is a growing interest in understanding the theoretical foundations that underlie its effectiveness. To gain theoretical insights into deep learning, numerous studies have been conducted from various perspectives, including approximation theory and statistical learning theory. The statistical learning community has extensively investigated the generalization properties of Deep Neural Networks (DNNs), \citep{Barron94estimation,Bartlett2017MarginBoundsNNs,Srebro2018PACMarginBoundsNNs,ZhangUnderstandingDL2017,SchmidtHieberDNN,Suzuki18DNNkerels,Imaizumi19DNN,suzuki2019adaptivity}. Specifically, \cite{SchmidtHieberDNN} and \cite{suzuki2019adaptivity} demonstrated that estimators based on sparsely connected DNNs with ReLU activation functions and carefully selected architectures can achieve minimax estimation rates (up to logarithmic factors) under classical smoothness assumptions on the regression function. Additionally, \cite{Imaizumi19DNN} and \cite{Suzuki2019Superiority} highlighted the superiority of DNNs over linear operators in certain scenarios, where DNNs achieve the minimax rate of convergence while alternative methods fall short.

On the Bayesian front, theoretical research has been relatively limited. Studies such as \cite{Rockova2018,Suzuki18DNNkerels,lee2022asymptotic,kong2023masked} have explored the  properties of the posterior distribution, while \cite{Vladimirova2019PriorsBNNsUnits} investigated the regularization effects of prior distributions at the unit level. Furthermore, \cite{cherief2020convergence,bai2020efficient} examined properties related to variational inference.
 
In the realm of classification, theoretical investigations of deep learning classifiers have been undertaken quite recently in various setups. Such as \cite{kim2021fast,shen2022approximation,hu2022minimax,bos2022convergence,wang2022minimax,kohler2022rate,meyer2023optimal,kohler2020statistical}. While a variational Bayes classification for dense deep neural networks has been studied \cite{bhattacharya2024comprehensive}, however, the prediction error for classification using Bayesian sparse deep learning remains unexplored. This study aims to address this gap by presenting an effort to shed light on this aspect.

We investigate a probabilistic approach for a sparse deep neural network classifier, where the network's architecture is governed by a Gibbs posterior incorporating a risk concept based on the hinge loss \citep{zhang2004statistical}. The hinge loss is commonly utilized for classification tasks in deep learning \citep{kim2021fast,shen2022approximation,hu2022minimax,bos2022convergence,wang2022minimax,kohler2022rate,kohler2020statistical}. Our analysis of prediction errors relies on the PAC-Bayes bound technique, originally introduced by \cite{McA,STW}. This technique aims to provide numerical generalization certificates for various probabilistic machine learning algorithms, and subsequently, \cite{catoni2004statistical,catonibook} demonstrated its utility in deriving oracle inequalities and rates of convergence for different methods. This methodology shares significant connections with the ``information bounds" presented by \cite{zha2006,russo2019much}. For a comprehensive exploration of this topic, we recommend consulting \cite{guedj2019primer,alquier2021user}.

Utilizing Spike-and-Slab priors, our main results are presented in terms of non-asymptotic oracle inequalities. These inequalities demonstrate that the prediction error of our proposed method is as good as the best possible errors. Our results are applicable in a general setting. Specifically, we show that by considering different architectures, our results can achieve minimax optimal rates for H\"older smooth functions in both low and high dimensions, up to a logarithmic factor. This alignment with works in the frequentist literature, such as \cite{kim2021fast} and \cite{wang2022minimax}, underscores the robustness of our findings. Moreover, our inclusion of an additional logarithmic term yields slight improvements over previous works.
Furthermore, we also propose and analyze an automated model selection approach aimed at optimally choosing a network architecture with guaranteed optimality.

The rest of the paper is structured as follows: Section \ref{sc_problem_method} introduces the problem, describes the deep neural network, and presents our proposed method. Our main results are outlined in Section \ref{sc_main_results}. All technical proofs are provided in Section \ref{sc_proofs}.

\subsubsection*{Notations}
Let us introduce the notations and the statistical framework we adopt in this paper. For any vector $x=(x_1,...,x_d) \in [-1,1]^d$ and any real-valued function $f$ defined on $[-1,1]^d$, $d>0$, we denote:
$
\|x\|_\infty = \max_{1\leq i \leq d} |x_i|
\textnormal{,}
\,
\|f\|_2 = ( \int f^2)^{1/2}
\textnormal{ and }
\,
\|f\|_\infty = \sup_{y\in [-1,1]^d} |f(y)|.
$ Put $ (a)_+ := \max(a,0),\forall a \in \mathbb{R} $.
For any $\mathbf{k}\in \{0,1,2,...\}^d$, we define $|\mathbf{k}|=\sum_{i=1}^d k_i$ and the mixed partial derivatives when all partial derivatives up to order $|\mathbf{k}|$ exist:
$
D^{\mathbf{k}} f(x) = \frac{\partial^{|\mathbf{k}|}f}{\partial^{k_1}x_1...\partial^{k_d}x_d} (x) .
$
For $\beta>0$, let  $\lfloor \beta \rfloor $ denote the largest integer strictly smaller than $\beta$ and . Then $f$ is said to be $\beta$-H\"older continuous \citep{tsybakov2008} if all partial derivatives up to order $\lfloor \beta \rfloor$ exist and are bounded. The norm $ \|f\|_{\mathcal{C}_\beta} $ of the H\"older space $\mathcal{C}_\beta = \{f : \|f\|_{\mathcal{C}_\beta}<+\infty\}$ is defined as:
\\
$
\|f\|_{\mathcal{C}_\beta} := \max_{|\mathbf{k}|\leq \lfloor \beta \rfloor} \|D^{\mathbf{k}} f\|_\infty + \max_{|\mathbf{k}| = \lfloor \beta \rfloor} \sup_{x,y \in [-1,1]^d, x\ne y} \frac{|D^\mathbf{k} f(x)-D^\mathbf{k} f(y)|}{\|x-y\|_\infty^{\beta-\lfloor \beta \rfloor}}.
$

\section{Problem and method}
\label{sc_problem_method}
\subsection{Problem statement}
\label{sc_modelstatement}

This study examines a general binary classification problem. Given a high-dimensional feature vector \( x \in \mathbb{R}^d \), the class label outcome \( Y \) is a binary variable derived from \( \{-1,1\} \) with an unknown probability \( p(x) \). Specifically, \( Y \) is equal to $ 1 $ with a probability of \( p(x) \), and $ -1 $ with a probability of \( 1 - p(x) \) when conditioned on \( x \). The accuracy of a classifier \( \eta \) is determined by the prediction or misclassification error, defined as
\[
R(\eta) = \mathbb{P}(Y \neq \eta(x)).
\]
However, the probability function $p( x)$ is unknown and the resulting classifier $\hat{\eta}( x)$ should be designed from the data $D_n $: a random sample of $n$ independent observations $( x_1,Y_1),\ldots, ( x_n,Y_n) $. 
The design points $ x_i$ may be considered as fixed or random.  
The corresponding (conditional) prediction error of $\hat{\eta}$ is 
$
R(\hat{\eta})
=
\mathbb{P} (Y \neq \hat{\eta}( x) \, |D_n)
$.
In this study, our objective is to construct a classifier \( \hat{\eta} \) with a prediction error that closely approximates the ideal Bayes error \( R^*:= \inf R(\eta) \), \citep{vapnik,devroye1997probabilistic}.

\subsection{Deep neural networks}

We define a deep neural network as any function \( f_\theta: \mathbb{R}^d \rightarrow \mathbb{R} \), recursively defined as:
\\
$
\begin{cases}
x^{(0)}:=x, 
\\
x^{(\ell)}:= \rho(A_\ell x^{(\ell-1)} + b_\ell) \hspace{0.3cm} \text{for} \hspace{0.1cm} \ell=1,...,L-1,
 \\
f_\theta(x):=A_L x^{(L-1)} + b_L
\end{cases}
$
\\
where $L\geq 3$ and $\rho$ is an activation function acting componentwise. 

For example, we can employ the ReLU activation function, \( \rho(u) = \max(u, 0) \). Each \( A_\ell \in \mathbb{R}^{D_\ell \times D_{\ell-1}} \) serves as a weight matrix. Its entry at the \( (i,j) \)-th position, called edge weight, connects the \( j \)-th neuron from the \( (\ell-1) \)-th layer to the \( i \)-th neuron of the \( \ell \)-th layer. And, each node vector \( b_\ell \in \mathbb{R}^{D_\ell} \) acts as a shift vector. Its \( i \)-th entry signifies the weight linked with the \( i \)-th node of the \( \ell \)-th layer.

Put $D_0=d$ the number of units in the input layer, $D_{L}=1$ the number of units in the output layer and $D_\ell=D$ the number of units in the hidden layers. The architecture of the network is characterized by its number of edges $S$ (the total number of nonzero entries in matrices $A_\ell$ and vectors $b_\ell$), its number of layers $L\geq 3$ (excluding the input layer), and its width $D \geq 1$. We have $ S\leq T := \sum_{\ell=1}^L D_\ell(D_{\ell-1}+1) $, the total number of coefficients in a fully connected network. 

At this point, we assume that \( S \), \( L \), and \( D \) are constants. We require \( d \leq D \) and \( T \leq LD(D+1) \). It is assumed that the absolute values of all coefficients are bounded above by a constant \( C_B \geq 2 \). The parameter for a Deep Neural Network (DNN) is denoted as \( \theta = \{ (A_1, b_1), \ldots, (A_L, b_L) \} \), and we define \( \Theta_{S,L,D} \) as the set of all possible parameters. Alternatively, we will also consider the stacked coefficients parameter, \( \theta = (\theta_1, \ldots, \theta_T) \).

\subsection{EWA procedure using empirical hinge loss}

We construct a classifier $ \hat{\eta}( x): = {\rm sign} (f_\theta (x)) $.

We adopt a PAC-Bayesian approach \citep{alquier2021user} and consider an exponentially weighted aggregate (EWA) procedure. Let's consider the following posterior distribution:
\begin{align}
\label{eq_mainporsterior}
\hat{\rho}_\lambda(\theta)
\propto
\exp[-\lambda r^h_n(\theta)] \pi(\theta)
\end{align}
where $\lambda>0$ is a tuning parameter that will be discussed later and $\pi(\theta)$ is a prior distribution, given in \eqref{eq_the_prior}, that promotes (approximately) sparsity on the parameter vector $ \theta $. In this paper, we primarily focus on the hinge loss resulting in the following hinge empirical risk:
\begin{align*}
r^h_n (\theta) 
= 
\frac{1}{n}\sum_{i=1}^n ( 1 - Y_{i} \,
f_\theta ( x_i) )_+ \, 
.
\end{align*}
It is noteworthy that several studies in deep learning for classification have employed hinge loss as a surrogate loss to facilitate computation, as in \citep{kim2021fast,wang2022minimax}.

The EWA procedure has found application in various contexts in prior works \citep{dalalyan2018exponentially,dalalyan2012sparse,dalalyan2020exponential,dalalyan2008aggregation}. The term $ \hat{\rho}_\lambda $ is also referred to as the Gibbs posterior \citep{AlquierRidgway2015,catonibook}. The incorporation of $ \hat{\rho}_\lambda $ is driven by the minimization problem presented in Lemma \ref{lemma_donvara}, rather than strictly adhering to conventional Bayesian principles.

\subsubsection{Sparsity prior}

We adopt a spike-and-slab prior \(\pi\) \citep{Castillo2015SS,cherief2020convergence}  over the parameter space \( \Theta_{S,L,D} \). The spike-and-slab prior is recognized as a relevant alternative to dropout in Bayesian deep learning, as discussed in \cite{Rockova2018}. More specifically, the prior is defined hierarchically as follows. Initially, a binary vector of indicators \( \gamma = (\gamma_1, \ldots, \gamma_T) \) is sampled from the set \( \mathcal{S}^S_T \) of \( T \)-dimensional binary vectors with precisely \( S \) non-zero entries. Subsequently, given \( \gamma_t \) for each \( t = 1, \ldots, T \), we apply a spike-and-slab prior to \( \theta_t \), which returns \( 0 \) if \( \gamma_t = 0 \) and a random sample from a uniform distribution on \( [-C_B, C_B] \) otherwise:
\begin{eqnarray}
\label{eq_the_prior}
\begin{cases}
\gamma \sim \mathcal{U}(\mathcal{S}^S_T), \\
\theta_t|\gamma_t \sim \gamma_t \hspace{0.1cm} \mathcal{U}([ -C_B , C_B ]) + (1-\gamma_t)\delta_{\{0\}}, \hspace{0.2cm} t=1,...,T
\end{cases}
\end{eqnarray}
where $\delta_{\{0\}}$ is a point mass at $0$ and $\mathcal{U}([ -C_B , C_B ])$ is a uniform distribution on $[ -C_B , C_B ]$. We recall that the sparsity level $S$ is fixed here and that this assumption will be relaxed in Section \ref{sc_modelselection}.

\begin{remark}
We opt for uniform distributions for simplicity, as done in related studies \citep{Rockova2018,Suzuki18DNNkerels}. However, Gaussian distributions can also be employed when dealing with an unbounded parameter set \( \Theta_{S,L,D} \), as indicated in \cite{cherief2020convergence}.
\end{remark}

\section{Main results}
\label{sc_main_results}

\subsection{Assumptions}
The result in this section applies to a wide range of activation functions, including the popular ReLU activation and the identity map. In the following, we make the following assumption regarding the activation function.

\begin{assume}
	\label{asm1}
Assume that the activation function $\rho$ is $1$-Lispchitz continuous with respect to the absolute value: $|\rho(x)|\leq|x|$ for any $x\in\mathbb{R}$.
\end{assume}

Put $ R^* := \min_\theta\mathbb{P}(Y \neq {\rm sign} (f_\theta (x)) ) \, $; $ \theta^*:= \arg\min_\theta\mathbb{P}(Y \neq {\rm sign} (f_\theta (x)) )\, $; $ r_n(\theta) 
= n^{-1} \sum_{i=1}^n \mathbf{1} \{ Y_{i} f_\theta (x_i) < 0 \} \, $, and $ r^*_n := r_n(\theta^*) $.
\begin{assume}
	\label{assume_bound_on_thetruebayes}
	We assume that there is a constant $ C' >0 $ such that $r^h_n(\theta^*) \leq (1+C')r^*_n $.
\end{assume}

Assumption \ref{assume_bound_on_thetruebayes} implies that for the underlying parameter \( \theta^* \), the function \( f_{\theta^*} (x) \) is bounded above by a constant.

\begin{assume}[Low-noise assumption] 
	\label{assume_margin}
	We assume that there is a constant $C \geq 1 $ such that: 
	$
	\mathbb{E}\left[\left(
	\mathbf{1}_{Y f_{\theta}(x) \leq 0} - 
	\mathbf{1}_{Y f_{\theta^*}(x) \leq 0} \right)^2\right] 
	\leq 
	C[R(\theta)-R^* ].
	$
\end{assume}

Assumption \ref{assume_margin} serves as the noise condition, capturing the difference between a classifier's performance and random guessing. This assumption is also commonly referred to as the margin assumption, as discussed in \cite{mammen1999smooth}, and further explored in \cite{tsybakov2004optimal,bartlett2006convexity}.

\subsection{General results}

Here, we present general results that do not require $ f_{\theta} $ to be $\beta$-H\"older and we consider any structure $(S,L,D)$.

\begin{theorem}
	\label{thm_slowrate}
Given Assumption \ref{asm1} and Assumption \ref{assume_bound_on_thetruebayes}, for \( \lambda = \sqrt{n} \), we find that with a probability of at least \( 1-2\epsilon \), where \( \epsilon \in (0,1) \), the following holds:
		\begin{align*}
	\int R d\hat{\rho}_\lambda 
	\leq 
	(1+2C')R^*  	
	+
c	\frac{	S \log \left(  \frac{nT D^L [ (d+1)L + 1] }{ S (D-1)}
	\right)
	+\log(1/\epsilon)  }{ \sqrt{n } } 
	,
	\end{align*}
	where $ c $ depends only ong $ C_B, C' $.
\end{theorem}

In addition to the results in Theorem \ref{thm_slowrate} for the integrated error, we can derive a result for, $ \hat{\theta} \sim \hat{\rho}_{\lambda} $, a stochastic classifier sampled from the Gibbs-posterior \eqref{eq_mainporsterior}. 

\begin{theorem}
	\label{thrm_contraction_slow} 
	Under the same assumptions for Theorem~\ref{thm_slowrate}, and the same definition for $\lambda $, for any small $\varepsilon \in (0,1) $. one has that 
	$ 
	\mathbb{E} \Bigl[ \mathbb{P}_{ \hat{\theta} \sim \hat{\rho}_{\lambda}} 
	( 	\hat{\theta} \in \mathcal{B}	 ) \Bigr] 
	\geq 
	1- \varepsilon
	,
	$
	where 
	\\
$
	\mathcal{B}	
	= 
	\Biggl\{ \theta \in \Theta_{S,L,D} : 
	R 
	\leq  
	(1+2C')R^* 	+
	c	\frac{	S \log \left( \frac{ T n D^L [ (d+1) L  + 1  ] }{ S( D-1) }
		 \right)
		+\log(2/\epsilon)  }{ \sqrt{n } }  
	\Biggl\}
	.
$
\end{theorem}

The oracle inequality presented in Theorem \ref{thm_slowrate} links the integrated prediction risk of our method to the minimum attainable risk. By incorporating additional assumptions, the bounds provided in Theorem \ref{thm_slowrate} can be improved.

\begin{theorem}
	\label{thm_fastrate}
Suppose  Assumption \ref{asm1}, Assumption \ref{assume_bound_on_thetruebayes} and Assumption \ref{assume_margin} are satisfied. For $ \lambda = 2 n/(3C + 2) $, we find that with a probability of at least \( 1-2\epsilon \), where \( \epsilon \in (0,1) \), the following holds:
	\begin{align*}
\int R d\hat{\rho}_\lambda  
\leq 	
(1+3C' ) R^*  	
+ 
c'\frac{ S \log \left( \frac{ T n D^L [ (d+1) L  + 1  ] }{ S( D-1) }
	\right) + \log\left(1/\epsilon\right)
}{n} 
.
\end{align*}
	where $ c' $ is a universal constant depending only on $ C,C',C_B $.
\end{theorem}

In contrast to Theorem \ref{thm_slowrate}, the bound in Theorem \ref{thm_fastrate} exhibits a faster rate, scaling as $1/n$ rather than $1/\sqrt{n}$. These bounds enable a comparison between the out-of-sample error of our method and the optimal error $R^*$.

For the noiseless case, i.e. $ Y = {\rm sign} ( f_{\theta^*}(x) ) $ almost surely, then $ R^* =0 $ and consequently that
$
\mathbb{E}\left[\left(\mathbf{1}_{Yf_{\theta}(x)\leq 0} 
- 
\mathbf{1}_{Y f_{\theta^*}(x) \leq 0} \right)^2\right]
= 
\mathbb{E}\left[\mathbf{1}_{Yf_{\theta}(x) \leq 0}^2\right]
= 
\mathbb{E}\left[\mathbf{1}_{Yf_{\theta}(x) \leq 0}\right] 
=  
R(\theta)-R^* .
$
Thus, Assumption \ref{assume_margin} is satisfied with $C=1 $.
\begin{corollary}
	\label{cor_noiseless}
 In the noiseless case, for $ \lambda = 2 n/5 $, with probability at least \( 1-2\epsilon, \epsilon \in (0,1) \), Theorem \ref{thm_fastrate} leads to 
	\begin{align*}
	\int R d\hat{\rho}_\lambda  
	\leq 	
	c'\frac{ S \log \left( \frac{ T n D^L [ (d+1) L  + 1  ] }{ S( D-1) }
		\right) + \log\left(1/\epsilon\right)
	}{n} 
	.
	\end{align*}
	where $ c' $ is a universal constant depending only on $ C',C_B $.
\end{corollary}

In analogy to Theorem \ref{thrm_contraction_slow}, we can establish that a stochastic classifier, $ \hat{\beta} \sim \hat{\rho}_{\lambda} $, drawn from our proposed pseudo-posterior in Equation \eqref{eq_mainporsterior} exhibits a fast rate. 
 
\begin{theorem}
	\label{thrm_contraction} 
	Under the same assumptions for Theorem~\ref{thm_fastrate}, and the same definition for \( \lambda \), along with any small \( \varepsilon \in (0,1) \), we find that
	$ 
	\mathbb{E} \Bigl[ \mathbb{P}_{ \hat{\theta} \sim \hat{\rho}_{\lambda}} 
	( 	\hat{\theta} \in \Omega 	) \Bigr] 
	\geq 
	1 - \varepsilon 
	,
	$
	where
	\\
$
	\Omega_n 	= 
	\Biggl\{ \theta \in \Theta_{S,L,D}: 
	R 
	\leq 
(1+3C' ) R^*  	
+ 
c'\frac{ S \log \left( \frac{ T n D^L [ (d+1) L  + 1  ] }{ S( D-1) }
	\right) + \log\left(2/\epsilon\right)
}{n} 
	\Biggr\}
	.
$
\end{theorem}

In this section, we establish specific values for the tuning parameters $ \lambda $ in our proposed method within theoretical results for prediction errors. However, it is noted that these values may not necessarily be optimal for practical applications. Cross-validation can be employed in practical scenarios to fine-tune the parameters appropriately. Nevertheless, the theoretical values identified in our analysis offer insights into the scale of the tuning parameters when utilized in real-world situations.

\subsection{Rates in specific settings}

\subsubsection*{Low-dimension setup}

Now, let's explore the scenario of low-dimensional settings, where we assume that \( d \) is fixed and smaller than \( n \). We will derive prediction error rates for the noiseless case, considering the following architecture.
\begin{example}
	\label{ex_1}
	Let us assume that $f_{\theta^*} $ is $\beta$-H\"older smooth with $0<\beta<d$ and that the activation function is ReLU. We consider the architecture of \cite{Rockova2018}, for some positive constant $C_D$ independent of $n$, that:
	$
	L
	= 8 + (\lfloor\log_2n\rfloor + 5)(1 + \lceil\log_2 d\rceil),
	D 
	= C_D \lfloor n^{\frac{d}{2\beta+d}}/\log n \rfloor,
	S 
	\leq 94 d^2 (\beta+1)^{2d} D (L + \lceil\log_2 d\rceil).
	$
\end{example}

\begin{proposition}
	\label{propos_lowdim}
	For the architecture as in Example \ref{ex_1}, then the rate in Corollary \ref{cor_noiseless} is of order $	
	 n^{\frac{-2\beta}{2\beta+d}}  \log n $.
\end{proposition}

\begin{remark}
The rate provided in Proposition \ref{propos_lowdim}, adjusted by a \( \log n \) factor, is minimax optimal \cite{tsybakov2004optimal}. Additionally, it is noteworthy that our \( \log n \) factor represents a slight enhancement over the \( \log^3 n \) factor achieved by the authors in \cite{kim2021fast} for obtaining the minimax optimal rate.
\end{remark}

\subsubsection*{High-dimensional setup}

Next, we delve into the realm of high-dimensional settings, where we assume that \( d > n \). We will consider the following architecture.

\begin{example}
	\label{ex_2}
	Let us assume that $f_{\theta^*} $ is $\beta$-H\"older smooth with $0<\beta<d$ and that the activation function is ReLU. We consider the architecture of \cite{wang2022minimax}:
	$
	L
	\asymp \log n,
	D 
	\asymp d,
	S 
\asymp n^{\frac{d}{2\beta+d}} \log n
.
	$
\end{example}
\begin{proposition}
	\label{propos_high_dim}
	For the architecture as in Example \ref{ex_2}, then the rate in Corollary \ref{cor_noiseless} is of order
$	
	 n^{\frac{-2\beta}{2\beta+d}} \log n \log d
$.
\end{proposition}

\begin{remark}
	The rate, $  n^{\frac{-2\beta}{2\beta+d}} \log d $, provided in Proposition \ref{propos_high_dim}, up to a \( \log n \) factor, is minimax optimal in a high dimensional setting as in \cite{wang2022minimax}. Additionally, it is noteworthy that our \( \log n \) factor represents a slight enhancement over the \(  \log^2 n \log d \) factor achieved by the authors in \cite{wang2022minimax} for obtaining the minimax optimal rate.
\end{remark}

\subsection{Architecture design via model selection}
\label{sc_modelselection}

Let  $\mathcal{M}_{S,L,D}$ denote the model associated with the parameter set $\Theta_{S,L,D}$. We consider a countable number of models, and introduce prior beliefs $ {\rm p}_{S,L,D} $ over the sparsity, the depth and the width of the network, hierarchically.

More specifically, $ {\rm p}_{S,L,D} = {\rm p}_{S|L,D} {\rm p}_{D|L} {\rm p}_L $ where $ {\rm p}_L=2^{-L}$, $ {\rm p}_{D|L} $ follows a uniform distribution over $\{d,...,\max(e^L,d)\}$ given $L$, and $ {\rm p}_{S|L,D} $ is a uniform distribution over $\{1,...,T\}$ given $L$ and $D$. Here, 
$ T $ represents the number of coefficients in a fully connected network. This particular choice is sensible as it allows to consider any number of hidden layers and (at most) an exponentially large width with respect to the depth of the network. We continue to employ spike-and-slab priors in \eqref{eq_the_prior} on $ \theta_{S,L,D} \in \Theta_{S,L,D} $ given model $\mathcal{M}_{S,L,D} $,  now denoted as $ \pi_{(S,L,D)} (\theta) $.

The Gibbs posterior is now as
$
\hat{\rho}_\lambda^{(S,L,D)} (\theta)
\propto
\exp[-\lambda r^h_n(\theta)] \pi_{(S,L,D)}(\theta)
$.
We formulate the following optimization problem to select the architecture of the network.
\begin{align}
\label{eq_optimal_modelselection}
(\hat{S}, \hat{L}, \hat{D})
=
\arg \min_{S,L,D}
\left[ \int 
r^h_n {\rm d}\hat{\rho}_\lambda^{(S,L,D)}  
+ 
\frac{\textnormal{KL}(\hat{\rho}_\lambda^{(S,L,D)},\pi_{(S,L,D)}) + \log(1/{\rm p}_{S,L,D}) }{\lambda}  
\right] 
\end{align}
The following theorem demonstrates that this strategy for selecting the architecture yields identical results to those in our main theorem, Theorem \ref{thm_fastrate}.

\begin{theorem}
	\label{thm_fasrate_modelselect}
	Suppose  Assumption \ref{asm1}, Assumption \ref{assume_bound_on_thetruebayes} and Assumption \ref{assume_margin} are satisfied. For $ \lambda = 2 n/(3C + 2) $, we find that with a probability of at least \( 1-2\epsilon \), where \( \epsilon \in (0,1) \), the following holds:
	\begin{align*}
	\int R d\hat{\rho}_\lambda^{(\hat{S}, \hat{L}, \hat{D})}  
	\leq 	
	(1+3C' ) R^*  	
	+ 
\inf_{S,L,D}
\left[
	c'\frac{ S \log \left( \frac{ T n D^L [ (d+1) L  + 1  ] }{ S( D-1) }\right) 
+ 
		\log(\frac{1}{{\rm p}_{S,L,D}}) 
	+ \log(\frac{1}{\epsilon})
	}{n} 
\right]
	.
	\end{align*}
	where $ c' $ is a universal constant depending only on $ C,C',C_B $.
\end{theorem}

Theorem \ref{thm_fasrate_modelselect} illustrates that once the complexity term \( \log(1/{\rm p}_{S,L,D})/n \), which reflects the prior beliefs, falls below the effective error \( r_n^{S,L,D} := S \max \left( \log(L), L\log(D), \log (nd) \right)/n \), the selection method in \eqref{eq_optimal_modelselection} adaptively achieves the best possible rate. This scenario leads to (near-)minimax rates for H\"older smooth functions and selects the optimal architecture. Notably, for the prior beliefs choice \(\pi_L=2^{-L}\), \(\pi_{D|L}=1/(\max(e^L,d)-d+1)\), \(\pi_{S|L,D}=1/T\), we obtain, up to some absolute constant \( K>0 \), that:
\[
\frac{\log(\frac{1}{{\rm p}_{S,L,D}})}{n} 
\leq 
K \frac{ \max(\log (D), \log L , L, \log (d)) }{n}
\]
which is lower than \( r_n^{S,L,D} \) (up to a factor).

It is worth noting that various model selection approaches for choosing the architecture of the network have been proposed in the study of deep learning, such as \cite{kim2021fast,cherief2020convergence}.

\section{Proof}
\label{sc_proofs}

\subsection*{Proof for slow rate}

\begin{proof}[\bf Proof of Theorem~\ref{thm_slowrate}]
	\text{}
	\\
	\textit{Step 1:}
	
	Put
	$
	U_{i} 
	=  
	\mathbf{1}_{Y_{i} f_{\theta}(x_i) \leq 0} - 
	\mathbf{1}_{Y_{i} f_{\theta^*}(x_i) \leq 0} 
	,
	$
	then, $ -1 \leq U_i \leq 1 $. One can utilize Hoeffding's Lemma \ref{lemma_hoeffding} to derive that
	\begin{align*}
	\mathbb{E}\exp\left\{\lambda[R( \theta )-R^* ]-\lambda[
	r_n( \theta )-r_n^* ] 
	\right\}
	\leq
	\exp\left( \lambda^2 / 2n  
	\right)
	.
	\end{align*}
Integrating with respect to $ \pi $
	and then applying Fubini's theorem, one gets, for any $\lambda\in(0,n) $, that
	\begin{align}
	\label{eq_t_contraction_slow}
	\mathbb{E} \int \exp\left\{\lambda[R( \theta )-R^* ]-\lambda[
	r_n( \theta )-r_n^* ] 
	-
	\frac{\lambda^2 }{2n} 
	\right\} 
	{\rm d}\pi( \theta ) 
	\leq 
	1
	,
	\end{align}	
	Consequently, using Lemma \ref{lemma_donvara},
	\begin{align*}
	\mathbb{E} \exp \left\lbrace \sup_{\rho} \int 
	\left\{ 
	\lambda [R( \theta )-R^* ] 
	-
	\lambda[r_n( \theta )-r_n^* ]
	-
	\frac{\lambda^2 }{2n} 
	\right\}
	\rho (d \theta) - 
	\textnormal{KL}(\rho,\pi) \right\rbrace 
	\leq 1.
	\end{align*}
Applying Markov's inequality, for \( \epsilon \in (0,1) \), we obtain that:
	\begin{align*}
	\mathbb{P} \left\{ \sup_{\rho} \int 
	\left[ 
	\lambda 
	[R( \theta )-R^* ]-\lambda [r_n( \theta )-r_n^* ] 
	-
	\frac{\lambda^2 }{2n}  \right]
	\rho (d\theta ) 
	- \textnormal{KL}(\rho,\pi) +\log \epsilon>0 \right\}  
\leq 
\epsilon.
	\end{align*}
Subsequently, considering the complement, we derive that with a probability of at least \(1-\epsilon\), the following holds:
	\begin{align*}
	\forall \rho, \quad \lambda \int  
	[R( \theta )-R^* ]\rho (d\theta ) 
	\leq 
	\lambda \int 
	[r_n( \theta )-r_n^* ]  \rho (d\theta ) + \textnormal{KL}(\rho,\pi) 
	+
	\frac{\lambda^2 }{2n} 
	+ 
	\log 
	\frac{1}{\epsilon}.
	\end{align*}
	Now, note that as $ r^h_n \geq r_n $ and
	as it stands for all $\rho$ then the right hand side can be minimized and, from Lemma \ref{lemma_donvara}, the minimizer over $\mathcal{P}(\Theta_{S,L,D}) $ is $\hat{\rho}_\lambda$.	Thus we get, when $\lambda > 0 $,
	\begin{align}
	\label{eq_pacBaeys_slowrate}
	\int R d\hat{\rho}_\lambda 
	\leq 
	R^*  + 
	\inf_{\rho \in \mathcal{P}(\Theta_{S,L,D}) } 
	\left[ \int 
	r^h_nd\rho + \frac{1}{\lambda} \textnormal{KL}(\rho,\pi) 
	\right] 
	- 
	r_n^*  
	+
	\frac{\lambda }{2n} 
	+ 
	\frac{1}{\lambda} 
	\log\frac{1}{\epsilon}
	.
	\end{align}	
\noindent	
Step 2:
	\\ 
	First, we have that,
	\begin{align}
	\int r^h_n (\theta) \rho (d\theta) 
	&	=  	
	\frac{1}{n} \int \sum_{i=1}^n ( 1 - Y_{i} f_{\theta}(x_i) )_{_+} \,  \rho (d\theta) 
	\nonumber
	\\
	& \leq
	\frac{1}{n} \left[ \sum_{i=1}^n 
	( 1 -  Y_{i} f_{\theta^*}(x_i) )_+
	+
	\int \sum_{i=1}^n  
	\left| 
	 f_{\theta}(x_i) - f_{\theta^*}(x_i) 
	\right|  
	\rho (d\theta)  \right]
	\nonumber
	\\
	& \leq
	r^h_n(\theta^{*}) 
	+ 	
	\frac{1}{n} \sum_{i=1}^n 	\int 
	\left| 
f_{\theta}(x_i) - f_{\theta^*}(x_i) 
\right|
	\rho (d\theta)
	.
	\label{eq_boundforhinge}
	\end{align}
	And using Lemma \ref{lem_boundfor_hinge},
	\begin{align*}
&	\int 
\left| 
f_{\theta}(x_i) - f_{\theta^*}(x_i) 
\right|
\rho (d\theta)
\\
&	\leq 
	\int 
r_L (\theta)
\rho (d\theta)
\\
& \leq
\int
	(C_B D )^{L-1}
\frac{ C_B D (d+1)-d}{ C_B D -1}
\sum_{u=1}^L \tilde{A}_u 
\rho (d\theta)
+ 
\int
\sum_{u=1}^L (C_B D )^{L-u} \tilde{b}_u 
\rho (d\theta)
\\
& \leq
(C_B D )^{L-1}
\frac{ C_B D (d+1)-d}{ C_B D -1}
\sum_{u=1}^L 
\int \tilde{A}_u 
\rho (d\theta)
+ 
\sum_{u=1}^L (C_B D )^{L-u}
\int \tilde{b}_u 
\rho (d\theta)
	. 
	\end{align*}
Here, we take \( \rho (d\theta) = q_n^*(\theta) \), as defined in equation \eqref{eq_special_qn}. We find that
$$
\int \tilde{A}_\ell q_n^*({\rm d} \theta) 
= 
\int \sup_{i,j} | A_{\ell,i,j}-A_{\ell,i,j}^* | q_n^*({\rm d} A_{\ell,i,j}) 
\leq 
s_n
,
$$
and 
$
\int \tilde{b}_u q_n^*({\rm d} \theta) 
= 
\int 
\sup_{j} |b_{u,j} - b_{u,j}^*|
 q_n^*({\rm d} b_{u,j} ) 
\leq 
s_n
.
$
Therefore,
	\begin{align}
	\int 
\left| 
f_{\theta}(x_i) - f_{\theta^*}(x_i) 
\right|
q_n^*({\rm d} \theta) 
& \leq
(C_B D )^{L-1}
\frac{ C_B D (d+1)-d}{ C_B D -1}
L s_n
+ 
s_n \sum_{u=1}^L (C_B D )^{L-u}
\nonumber
\\
& \leq
(C_B D )^{L-1}
\frac{ C_B D (d+1)-d}{ C_B D -1}
L s_n
+ 
s_n \sum_{\ell = 0}^{L-1} (C_B D )^{\ell}
\nonumber
\\
& \leq
(C_B D )^{L-1}
\frac{ C_B D (d+1)-d}{ C_B D -1}
L s_n
+ 
s_n \frac{(C_B D )^L -1}{ C_B D -1}
\nonumber
\\
& \leq
s_n
\left[
(C_B D )^{L}
\frac{ C_B D (d+1)-d}{(C_B D -1) C_B D} L 
+ 
 \frac{(C_B D )^L}{ C_B D -1}
\right]  \nonumber
\\
& \leq
s_n \frac{(C_B D )^L}{ C_B D -1}
\left[
(d+1) L 
+ 
1
\right]
\label{eq_bound_abs_term}
. 
\end{align}

	From Assumption \ref{assume_bound_on_thetruebayes}, we have   $r^h_n(\theta^*) \leq (1+C')r^*_n $. Plug in \eqref{eq_bound_abs_term} and \eqref{eq_lemma_on_KL},  into \eqref{eq_pacBaeys_slowrate}, we obtain
	\begin{multline*}
	\int R d\hat{\rho}_\lambda 
	\leq 
	R^* + 	C'r^*_n
	+
s_n \frac{(C_B D )^L}{ C_B D -1}
\left[ (d+1) L  + 1 \right]
\\
	+ 
	\frac{	S\log(TB) + \frac{S}{2} \log \left(\frac{1}{s_n^2}\right)  }{\lambda} 
	+
	\frac{\lambda }{2n} 
	+
	\frac{1}{\lambda} \log\left(\frac{1}{\epsilon}\right) 
	.
	\end{multline*}
	Then, we use Lemma \ref{lem_theo1}, with probability at least $ 1-2\epsilon $, to obtain that
	\begin{multline}
	\int R d\hat{\rho}_\lambda 
	\leq 
	(1+2C')R^*  
	+
	C'\frac{1}{n\varsigma }\log \frac{1}{\epsilon}
	+
s_n \frac{(C_B D )^L}{ C_B D -1}
\left[ (d+1) L  + 1 \right]
\\
+ 
\frac{	S\log(T C_B) + S \log \left(\frac{1}{s_n}\right)  }{\lambda} 
+
\frac{\lambda }{2n} 
+
\frac{1}{\lambda} \log\left(\frac{1}{\epsilon}\right) 
	\label{eq_used_forknwon_sstatr}
	.
	\end{multline}
	By taking  $ s_n= \frac{S}{n} \left[ \frac{(C_B D )^L}{ C_B D -1}
	\left[ (d+1) L  + 1 \right] \right]^{-1} $ and $ \lambda = \sqrt{n } $, we can obtain that
	\begin{multline*}
	\int R d\hat{\rho}_\lambda 
	\leq 
	(1+2C')R^*  	+
	C'	\frac{1}{n\varsigma }\log \frac{1}{\epsilon}
	+   \frac{S}{n}
	\\
	+ 
	\frac{	S \log \left( T C_B \frac{n}{S}  \frac{(C_B D )^L}{ C_B D -1}
		\left[ (d+1) L  + 1  \right] \right)  }{ \sqrt{n } } 
	+
	\frac{\log(1/\epsilon)}{\sqrt{n }} 
	.
	\end{multline*}
	The proof is completed.
\end{proof}

\begin{proof}[\textbf{Proof of Theorem \ref{thrm_contraction_slow}}]
	From \eqref{eq_t_contraction_slow}, we have that
	\begin{align*}
	\mathbb{E} \Biggl[ \int \exp \left\{ 
	\lambda[R( \theta )-R^* ] 
	-\lambda[r_n( \theta )-r_n^* ]
	- 
	\log \left[\frac{d\hat{\rho}_{\lambda}}{d \pi} (\theta)  \right]
	-
	\frac{\lambda^2 }{2n} 
	- 
	\log\frac{1}{\varepsilon}
	\right\}
	\hat{\rho}_{\lambda}(d \theta)
	\Biggr]
	\leq 
	\varepsilon
	.
	\end{align*}
We employ the Chernoff's trick, which is based on \( \exp(x) \geq 1_{\mathbb{R}_{+}}(x) \), leading to
	$
	\mathbb{E} \Bigl[ 
	\mathbb{P}_{\theta \sim \hat{\rho}_{\lambda}} 
	(\theta \in 	\mathcal{B}	 ) \Bigr]
	\geq 
	1- \varepsilon
	,
	$
	where
	$$
		\mathcal{B}	
	= 
	\left\{\theta : 	
	\lambda	[R( \theta )-R^* ] -\lambda[r_n( \theta )-r_n^* ]   
	\leq      
	\log \left[\frac{d\hat{\rho}_{\lambda}}{d \pi} (\theta)  \right]
	+ 	\frac{\lambda^2 }{2n} 
	+ \log\frac{2}{\varepsilon} \right\}.
	$$
Utilizing the definition of \( \hat{\rho}_\lambda \) and Lemma \ref{lemma_donvara}, and recognizing that \( r_n \leq r^h_n \), for \( \theta \in 	\mathcal{B}	 \), we find that:
	\begin{align*}
	\lambda	[ R(\theta) - R^* ]
	&\leq  
	\lambda\Bigl( r(\theta) - r_n^* \Bigr)  +       \log \left[\frac{d\hat{\rho}_{\lambda}}{d \pi} (\theta)  \right]
	+ 	\frac{\lambda^2 }{2n} 
	+ \log\frac{2}{\varepsilon}
	\\
	&    \leq  
	\lambda\Bigl( r^h_n(\theta) - r_n^* \Bigr)  +       \log \left[\frac{d\hat{\rho}_{\lambda}}{d \pi} (\theta)  \right]
	+ 	\frac{\lambda^2 }{2n} 
	+ \log\frac{2}{\varepsilon}
	\\
	& \leq 
	- \log\int\exp\left[
	-\lambda r^h_n (\theta) \right]\pi({\rm d} \theta) - \lambda r_n^*
	+ 	\frac{\lambda^2 }{2n} 
	+ \log\frac{2}{\varepsilon}
	\\
	& = \lambda\Bigl( \int r^h_n (\theta) \hat{\rho}_{\lambda}({\rm d} \theta) - r_n^* \Bigr)  +    \mathcal{K}(\hat{\rho}_\lambda,\pi)
	+ 	\frac{\lambda^2 }{2n} 
	+ \log\frac{2}{\varepsilon}
	\\
	& = 
	\inf_{\rho} \left\{ 
	\lambda\Bigl( \int r^h_n (\theta) \rho({\rm d}\theta) - r_n^* \Bigr)  +    \mathcal{K}(\rho,\pi)
	+ 	\frac{\lambda^2 }{2n} 
	+ \log\frac{2}{\varepsilon} \right\}.
	\end{align*}
	We upper-bound the right-hand side exactly as Step 2 in the proof of Theorem~\ref{thm_slowrate}. The result of the theorem is followed.	
\end{proof}

\subsection*{Proof for fast rate}

\begin{proof}[\bf Proof of Theorem~\ref{thm_fastrate}]
	\text{}
	\\
	\textit{Step 1:}
	
	Fix any $ \theta $ and put
	$
	U_{i} 
	=  
	\mathbf{1}_{Y_{i} f_{\theta}(x_i) \leq 0} - 
	\mathbf{1}_{Y_{i} f_{\theta^*}(x_i) \leq 0} 
	.
	$
	Under Assumption \ref{assume_margin}, we have that
	$
	\sum_{i} \mathbb{E}[U_{i}^{2}]
	\leq
	nC[R( \theta )-R^* ] .$
	Now, for any integer $k\geq 3$, as the 0-1 loss is bounded, it results in
	$$
	\sum_{i} \mathbb{E}\left[(U_{i})_+^{k}\right] \leq
	\sum_{i} \mathbb{E}\left[ |U_{i}|^{k-2} |U_{i}|^{2}\right]
	\leq
	\sum_{i} \mathbb{E}\left[ |U_{i}|^{2}\right].
	$$ 
	Thus, we can apply Lemma~\ref{lemmemassart} with  $v := 	nC[R( \theta )-R^* ] $, $w:=1 $ and $\zeta := \lambda/n $. We obtain, for any $\lambda\in(0,n) $,
	\begin{align*}
	\mathbb{E}\exp\{\lambda (
	[R( \theta )-R^* ]-[
	r_n( \theta )-r_n^* ] )
	\} 
	\leq 
	\exp \left\{ \frac{C\lambda^2[R( \theta )-R^* ] }{2n (1-\lambda/n)} 	 
	\right\},
	\end{align*}
By integrating with respect to \( \pi \) and subsequently applying Fubini's theorem, we find that	
\begin{align}
	\label{eq_to_contraction}
	\mathbb{E} \left[
	 \int \exp \left\{ 
	(\lambda- \frac{C\lambda^2 }{2n(1-\lambda/n)} )[R( \theta )-R^* ] 
	-\lambda[r_n( \theta )-r_n^* ]
	\right\} \pi(d\theta ) 
	\right]
	\leq 1
	.
	\end{align}
	Consequently, using Lemma \ref{lemma_donvara},
	\begin{align*}
	\mathbb{E} \left[
	 e^{ \sup_{\rho} \int \left\{ 
	(\lambda- \frac{C\lambda^2 }{2n(1-\lambda/n)} )[R( \theta )-R^* ] 
	-\lambda[r_n( \theta )-r_n^* ]
	\right\}
	\rho (d M) - 
	\textnormal{KL}(\rho,\pi)
}
\right]
	\leq 1.
	\end{align*}
By applying Markov's inequality and subsequently considering the complement, we establish that with a probability of at least \(1-\epsilon\), the following holds:
	\begin{align*}
	\forall \rho, \, 
	\left( \lambda-\frac{C\lambda^2 }{2n(1-\lambda/n)}\right) 
	 \int  
	[R( \theta )-R^* ]\rho (d\theta ) 
	\leq 
	\lambda \int 
	[r_n( \theta )-r_n^* ]  \rho (d\theta ) + \textnormal{KL}(\rho,\pi) + \log 
	\frac{1}{\epsilon}.
	\end{align*}
	Now, observe that $ r^h_n\geq r_n $,
	\begin{align*}
	\lambda \left[ \int r_n d\rho - r_n^* \right] + 
	\textnormal{KL}(\rho,\pi)
	\leq  
	\lambda \left[ \int r^h_nd\rho + \frac{1}{\lambda} 
	\textnormal{KL}(\rho,\pi) 
	\right] - \lambda r_n^* 
	.
	\end{align*}
	As it stands for all $\rho$ then the right hand side can be minimized and, from Lemma \ref{lemma_donvara}, the minimizer over $\mathcal{P}(\Theta_{S,L,D}) $ is $\hat{\rho}_\lambda$.	Thus we get, when $\lambda<2n/(C+2) $,
	\begin{align}
	\int R d\hat{\rho}_\lambda 
	\leq 
	R^*  + \frac{1}{1-\frac{C\lambda }{2n(1-\lambda/n)}} 
	\left\lbrace 
	\inf_{\rho \in \mathcal{P}(\Theta_{S,L,D}) } 
	\left[ \int 
	r^h_nd\rho + \frac{
		\textnormal{KL}(\rho,\pi)
	+ 	\log\frac{1}{\epsilon}  }{\lambda}  
	\right] - r_n^* 
	\right\rbrace
	\label{eq_PACBayes_fast}
	.
	\end{align}
	
	\textit{Step 2:}	
	\\
	From \eqref{eq_boundforhinge}, we have that,
	\begin{align*}
	\int r^h_n (\theta) \rho (d\theta) 
	\leq
	r^h_n(\theta^{*}) 
+ 	
\frac{1}{n} \sum_{i=1}^n 	\int 
\left| 
f_{\theta}(x_i) - f_{\theta^*}(x_i) 
\right|
\rho (d\theta)
	.
	\end{align*}
Here, we take \( \rho (\theta) = q_n^*(\theta) \), as defined in equation \eqref{eq_special_qn}. Then from \eqref{eq_bound_abs_term},
 and \eqref{eq_lemma_on_KL}, we deduce \eqref{eq_PACBayes_fast} into that,
	from Assumption \ref{assume_bound_on_thetruebayes}, as   $r^h_n(\theta^*) \leq (1+C')r^*_n $, we have that
	\begin{multline*}
	\int R d\hat{\rho}_\lambda 
	\leq 
	R^* + 
	\frac{1}{1-\frac{C\lambda }{2n(1-\lambda/n)}} 
	\left\lbrace 
	C' r_n^*
	+
s_n \frac{(C_B D )^L}{ C_B D -1}
\left[
(d+1) L + 1
\right]
\right.
\\
\left.
	+ 
	\frac{S\log(TC_B) + \frac{S}{2} \log \left(\frac{1}{s_n^2}\right) 
	+
\log\left(\frac{1}{\epsilon}\right)  }{\lambda} 
	\right\rbrace
	.
	\end{multline*}
	Taking $\lambda = 2 n/(3C + 2) $,
	we obtain:
	\begin{multline*}
		\int R d\hat{\rho}_\lambda  
	\leq 	
	R^*	+ 
	\frac{3}{2} 
	\left\lbrace 
C' r_n^*
+
s_n \frac{(C_B D )^L}{ C_B D -1}
\left[
(d+1) L + 1
\right]
\right.
\\
\left.
+ 
(3C+2) \frac{S\log(TC_B) + \frac{S}{2} \log \left(\frac{1}{s_n^2}\right) 
	+
	\log\left(\frac{1}{\epsilon}\right)  }{2n} 
\right\rbrace
	.
	\end{multline*}
	Then, we use Lemma~\ref{lem_theo1}, with probability at least $ 1-2\epsilon $, to obtain that
		\begin{multline*}
	\int R d\hat{\rho}_\lambda  
	\leq 	
	(1+3C' ) R^*  	+ 
	\frac{3}{2} 
	\left\lbrace 
	C' \frac{1}{n\varsigma }\log \frac{1}{\epsilon}
	+
	s_n \frac{(C_B D )^L}{ C_B D -1}
	\left[
	(d+1) L + 1
	\right]
	\right.
	\\
	\left.
	+ 
(3C+2)	\frac{S\log(TC_B) + S \log \left(\frac{1}{s_n} \right)
	+\log\left(1/\epsilon\right)
 }{ 2n } 
	\right\rbrace
	.
	\end{multline*}
	By setting  $ s_n= \frac{S}{n} \left[ \frac{(C_B D )^L}{ C_B D -1}
\left[ (d+1) L  + 1 \right] \right]^{-1} $, one finds that 
		\begin{multline*}
\int R d\hat{\rho}_\lambda  
\leq 	
(1+3C' ) R^*  	+ 
\frac{3}{2} 
\left\lbrace 
C' \frac{1}{n\varsigma }\log \frac{1}{\epsilon}
+
\frac{S}{n}
\right.
\\
\left.
+ 
(3C+2) \frac{ S \log \left( TC_B \frac{n}{S}  \frac{(C_B D )^L}{ C_B D -1}
	\left[ (d+1) L  + 1 \right] 
	\right) + \log\left(1/\epsilon\right)
}{2n} 
\right\rbrace
.
\end{multline*}
	The proof is completed.
\end{proof}

\begin{proof}[\textbf{Proof of Theorem \ref{thrm_contraction}}]
	From \eqref{eq_to_contraction}, we have that
	\begin{align*}
	\mathbb{E} \Biggl[  \int \exp \left\{ 
	(\lambda- \frac{C\lambda^2 }{2n(1-\lambda/n)} )[R( \theta )-R^* ] 
	-\lambda[r_n( \theta )-r_n^* ]
	- 
	\log \left[\frac{d\hat{\rho}_{\lambda}}{d \pi} (\theta)  \right]
	- \log\frac{1}{\varepsilon}
	\right\}
	\hat{\rho}_{\lambda}(d \theta)
	\Biggr]
	\leq 
	\varepsilon
	.
	\end{align*}
	Using Chernoff's trick, i.e. $\exp(x) \geq 1_{\mathbb{R}_{+}}(x)$, this gives:
	$
	\mathbb{E} \Bigl[ 
	\mathbb{P}_{\theta \sim \hat{\rho}_{\lambda}} 
	(\theta \in \Omega ) \Bigr]
	\geq 1- \varepsilon
	$
	where
	$$
	\Omega
	= 
	\left\{\theta : 	
	(\lambda- \frac{C\lambda^2 }{2n(1-\lambda/n)} )
	[R( \theta )-R^* ] 
	-\lambda[r_n( \theta )-r_n^* ]   
	\leq      
	\log \left[\frac{d\hat{\rho}_{\lambda}}{d \pi} (\theta)  \right]
	+ \log\frac{2}{\varepsilon} 
	\right\}
	.
	$$
	Using the definition of $\hat{\rho}_\lambda $ and noting that $ r_n \leq r^h_n $, for $ \theta \in \Omega $ we have
	\begin{align*}
	\Bigl(
	\lambda- \frac{C\lambda^2 }{2n(1-\lambda/n)} 
	\Bigr)
	\Bigl[ R(\theta) - R^* \Bigr]
	&\leq  
	\lambda\Bigl( r_n(\theta) - r_n^* \Bigr)  +       \log \left[\frac{d\hat{\rho}_{\lambda}}{d \pi} (\theta)  \right]
	+ \log\frac{2}{\varepsilon}
	\\
	&    \leq  
	\lambda\Bigl( r^h_n(\theta) - r_n^* \Bigr)  +       \log \left[\frac{d\hat{\rho}_{\lambda}}{d \pi} (\theta)  \right]
	+ \log\frac{2}{\varepsilon}
	\\
	& \leq 
	- \log\int\exp\left[
	-\lambda r^h_n (\theta) \right]\pi({\rm d} \theta) - \lambda r_n^*
	+ \log\frac{2}{\varepsilon}
	\\
	& = \lambda\Bigl( \int r^h_n (\theta) \hat{\rho}_{\lambda}({\rm d} \theta) - r_n^* \Bigr)  +    \mathcal{K}(\hat{\rho}_\lambda,\pi)
	+ \log\frac{2}{\varepsilon}
	\\
	& = 
	\inf_{\rho} \left\{ 
	\lambda\Bigl( \int r^h_n (\theta) \rho({\rm d}\theta) - r_n^* \Bigr)  +    \mathcal{K}(\rho,\pi)
	+ \log\frac{2}{\varepsilon} \right\}.
	\end{align*}
	We upper-bound the right-hand side exactly as Step 2 in the proof of Theorem~\ref{thm_fastrate}. The result of the theorem is followed.	
\end{proof}

\begin{proof}[Proof of Proposition \ref{propos_lowdim}]
	For certain constants \( K > 0 \), independent of \( n \), we have that
	\begin{align*}
	\frac{ S \log \left( \frac{ T n D^L [ (d+1) L  + 1  ] }{ S( D-1) }
		\right) }{n} 
	& 
	\leq 
	K \frac{ S \log \left( \frac{ LD(D+1) n D^L [ (d+1) L  + 1  ] }{ S( D-1) }
		\right) }{n} 
	\\
	& \leq 
	K \frac{ S \max  
		\left( \log(L), L\log(D), \log n
		\right) }{n} 
	\\
	& \leq 
	K \frac{n^{\frac{d}{2\beta+d}}}{n} \log n = 
	K n^{\frac{-2\beta}{2\beta+d}} \log n .
	\end{align*}
\end{proof}

\begin{proof}[Proof of Proposition \ref{propos_high_dim}]
	For certain constants \( K > 0 \), independent of \( n , d\), we have that
	\begin{align*}
	\frac{ S \log \left( \frac{ T n D^L [ (d+1) L  + 1  ] }{ S( D-1) }
		\right) }{n} 
	& 
	\leq 
	K \frac{ S \log \left( \frac{ LD(D+1) n D^L [ (d+1) L  + 1  ] }{ S( D-1) }
		\right) }{n} 
	\\
	& \leq 
	K \frac{ S \max  
		\left( \log(L), L\log(D), \log (nd)
		\right) }{n} 
	\\
	& \leq 
	K \frac{n^{\frac{d}{2\beta+d}} \log n \log d}{n} 
	= 
	K n^{\frac{-2\beta}{2\beta+d}} \log n \log d.
	\end{align*}
\end{proof}

\subsection*{Proof of Section \ref{sc_modelselection}}

\begin{proof}[Proof of Theorem \ref{thm_fasrate_modelselect}]
	From \eqref{eq_PACBayes_fast}, we have that
	\begin{align*}
	\int R d\hat{\rho}_\lambda^{(S,L,D)} 
	\leq 
	R^*  + \frac{1}{1-\frac{C\lambda }{2n(1-\lambda/n)}} 
	\left\lbrace 
	\inf_{\rho } 
	\left[ \int 
	r^h_nd\rho + \frac{
		\textnormal{KL}(\rho,\pi_{(S,L,D)})
		+ 	\log\frac{1}{\epsilon}  }{\lambda}  
	\right] - r_n^* 
	\right\rbrace
	.
	\end{align*}
	From Lemma \ref{lemma_donvara}, it holds that
	$ \hat{\rho}_\lambda^{(S,L,D)} $ is the minimizer of $  \int 
	r^h_nd\rho + \frac{
		\textnormal{KL}(\rho,\pi_{(S,L,D)})
	}{\lambda}  $, thus
	\begin{multline*}
	\int R d\hat{\rho}_\lambda^{(S,L,D)} 
	\leq 
	R^* + \frac{1}{1-\frac{C\lambda }{2n(1-\lambda/n)}} \times
	\\
	\left\lbrace 
	\int 
	r^h_nd \hat{\rho}_\lambda^{(S,L,D)} 
	+ 
	\frac{ \textnormal{KL}(\hat{\rho}_\lambda^{(S,L,D)} ,\pi_{(S,L,D)})
		+
		\log(\frac{1}{{\rm p}_{S,L,D}}) 
		+ 	\log\frac{1}{\epsilon}  }{\lambda}  
	- r_n^* 
	\right\rbrace
	.
	\end{multline*}
	Now using the definition of $ (\hat{S}, \hat{L}, \hat{D}) $ in \eqref{eq_optimal_modelselection}, we arrive at
	\begin{multline*}
	\int R d\hat{\rho}_\lambda^{(\hat{S}, \hat{L}, \hat{D})} 
	\leq 
	R^*  + \frac{1}{1-\frac{C\lambda }{2n(1-\lambda/n)}} \times
	\\
	\left\lbrace 
	\inf_{S,L,D}
	\left[ \int r^h_nd \hat{\rho}_\lambda^{(S,L,D)} 
	+
	\frac{\textnormal{KL}(\hat{\rho}_\lambda^{(S,L,D)},\pi_{(S,L,D)}) 
		+ 
		\log(\frac{1}{{\rm p}_{S,L,D}}) 
		+ 	\log\frac{1}{\epsilon} }{\lambda}  
	\right] - r_n^* 
	\right\rbrace
	,
	\end{multline*}
	using Lemma \ref{lemma_donvara}, we deduce that 
	\begin{multline*}
	\int R d\hat{\rho}_\lambda^{(\hat{S}, \hat{L}, \hat{D})} 
	\leq 
	R^*  + \frac{1}{1-\frac{C\lambda }{2n(1-\lambda/n)}} \times
	\\
	\left\lbrace 
	\inf_{S,L,D} \inf_{\rho } 
	\left[ \int r^h_nd \rho
	+
	\frac{\textnormal{KL}(\rho ,\pi_{(S,L,D)}) 
		+ 
		\log(\frac{1}{{\rm p}_{S,L,D}}) 
		+ 	\log\frac{1}{\epsilon} }{\lambda}  
	\right] - r_n^* 
	\right\rbrace
	.
	\end{multline*} 
	To obtain the result, proceed as in Step 2 in the proof of Theorem \ref{thm_fastrate}, page \pageref{eq_PACBayes_fast}. This completes the proof.
\end{proof}

\subsection*{Auxiliary lemmas}

\begin{lemma}
	\label{lem_theo1}[Lemma 6 in \cite{cottet20181}]
	For $\epsilon \in (0,1) $, with probability at least $1-\epsilon$, we have, for every $ \varsigma \in (0,1) $, that
$
	r_n^* 
	\leq 
	(1+\varsigma)	R^* +\frac{1}{n\varsigma }\log \frac{1}{\epsilon}
	,
$
	or we can have that $ r_n^* 
	\leq 
	2R^* +\frac{1}{n\varsigma }\log \frac{1}{\epsilon} $.
\end{lemma}

We will use the following version of the Bernstein's lemma, from \cite[page 24]{MR2319879}.
\begin{lemma}
	\label{lemmemassart} Let $U_{1}$, \ldots, $U_{n}$ be independent real
	valued random variables. Assuming that there exist two constants
	$v$ and $w$ such that
	$
	\sum_{i=1}^{n} \mathbb{E}[U_{i}^{2}] \leq v 
	$
	and that for all integers $k\geq 3$,
	$
	\sum_{i=1}^{n} \mathbb{E}\left[(U_{i})_+^{k}\right] \leq v k!w^{k-2}/2. 
	$
	Then, for any $ \zeta \in (0,1/w)$,
	$$ 
	\mathbb{E}
	\exp\left[\zeta \sum_{i=1}^{n}\left[U_{i}-\mathbb{E}U_{i}\right]
	\right]
	\leq 
	\exp\left(\frac{v\zeta^{2}}{2(1-w\zeta )} \right) .
	$$
\end{lemma}

\begin{definition}
	Put $q_n^*$ as
	\begin{align}
	\begin{cases}
	\gamma_t^* = \mathbb{I}(\theta^*_t\ne 0), 
	\\
	\theta_t \sim \gamma_t^* \hspace{0.1cm} \mathcal{U}([\theta_t^*-s_n,\theta^*_t+s_n]) + (1-\gamma_t^*)\delta_{\{0\}}, \, 
	\textnormal{ for each } \,
	t=1,...,T.
	\end{cases}
	\label{eq_special_qn}
	\end{align}
\end{definition}

\begin{lemma}
	\label{lem_boundon_KL}
	We have, for $ q_n^* $ in \eqref{eq_special_qn}, that
	\begin{align}
	\label{eq_lemma_on_KL}
	\textnormal{KL}(q_n^*\|\pi) 
	\leq 
	S\log(TC_B) + \frac{S}{2} \log\bigg(\frac{1}{s_n^2}\bigg) 
	.
	\end{align}
\end{lemma}
\begin{proof}
	The proof of this Lemma can be found in the proof of Theorem 2 in \cite{cherief2020convergence}, in particular in the third step of the proof.
\end{proof}

\begin{lemma}
	\label{lem_boundfor_hinge}
	Define the loss of the output layer as
	$
	r_\ell(\theta) := \sup_{x\in[-1,1]^d}  |f_\theta^L(x) - f_{\theta^*}^L(x)| = \sup_{x\in[-1,1]^d}  |f_\theta(x) - f_{\theta^*}(x)|.
	$
	We have, under Assumption \ref{asm1}, that for any $\ell=1,...,L$:
	\\
$
	r_\ell(\theta) 
	\leq 
	( C_B D)^{\ell -1}
	\frac{ C_B D(d+1)-d}{C_B D -1}
	\sum_{u=1}^\ell \tilde{A}_u 
	+ 
	\sum_{u=1}^\ell (C_B D )^{\ell-u} \tilde{b}_u 
	,
$
\\
	where $\tilde{A}_u=\sup_{i,j} |A_{u,i,j}-A_{u,i,j}^*|$ and 
	$
	\tilde{b}_u = \sup_{j} |b_{u,j} - b_{u,j}^*|
	$.
\end{lemma}

\begin{proof}
	The proof of this Lemma can be found in the proof of Theorem 2 in \cite{cherief2020convergence}, in particular in the second step of the proof.
\end{proof}

 We remind here a version of Hoeffding's inequality
for bounded random variables.
\begin{lemma}
	\label{lemma_hoeffding}
	Let $ U_i, i = 1,\ldots,n$ be $n$ independent random
	variables with $ a \leq U_i \leq b $ a.s., and
	$\mathbb{E}(U_i)=0$. Then, for any $\lambda>0$,
	$
	\mathbb{E}\exp\left(\frac{\lambda}{n}
	\sum_{i=1}^{n} U_i\right)
	\leq
	\exp\left(\frac{\lambda^2 (b-a)^2}{8n} \right) .$
\end{lemma}

\begin{lemma}[\cite{catonibook}{Lemma 1.1.3}]
	\label{lemma_donvara}
	
	Let $\mu \in\mathcal{P}(\Theta)$. For any measurable, bounded function $h:\Theta\rightarrow\mathbb{R}$ we have:
$
	\log \int {\rm e}^{h(\theta)} \mu({\rm d}\theta) = 
	\sup_{\rho\in\mathcal{P}(\Theta)}\left[\int h(\theta) \rho({\rm d}\theta) 
	-
	\textnormal{KL} (\rho , \mu)\right].
$
	Moreover, the supremum w.r.t $\rho$ in the right-hand side is
	reached for the Gibbs distribution,
	$
	\hat{\rho} (d\theta)
	\propto 
	\exp(h(\theta)) \pi (d\theta).
	$
\end{lemma}

\bmhead{Acknowledgements}
The author was supported by the Norwegian Research Council, grant number 309960, through the Centre for Geophysical Forecasting at NTNU.

\bmhead{Conflict of interest/Competing interests}
 The author declares no potential conflict of interests.





\end{document}